\documentstyle[12pt]{article}

\begin{document}

\newtheorem{Def}{Definition}[section]
\newtheorem{Th}{Theorem}[section]
\newtheorem{Prop}{Proposition}[section]
\newtheorem{Not}{Notation}[section]
\newtheorem{Lemma}{Lemma}[section]
\newtheorem{Rem}{Remark}[section]
\newtheorem{Cor}{Corollary}[section]

\def\s{\section}
\def\ss{\subsection}

\def\d{\begin{Def}}
\def\t{\begin{Th}}
\def\p{\begin{Prop}}
\def\n{\begin{Not}}
\def\la{\begin{Lemma}}
\def\r{\begin{Rem}}
\def\c{\begin{Cor}}
\def\ee{\begin{equation}}
\def\aa{\begin{eqnarray}}
\def\ya{\begin{eqnarray*}}
\def\bd{\begin{description}}

\def\ed{\end{Def}}
\def\et{\end{Th}}
\def\epo{\end{Prop}}
\def\en{\end{Not}}
\def\el{\end{Lemma}}
\def\er{\end{Rem}}
\def\ec{\end{Cor}}
\def\eee{\end{equation}}
\def\eaa{\end{eqnarray}}
\def\ey{\end{eqnarray*}}
\def\ebd{\end{description}}

\def\nn{\nonumber}
\def\bp{{\bf Proof.}\hspace{2mm}}
\def\qe{\hfill$\Box$}
\def\lj{\langle}
\def\rj{\rangle}
\def\dd{\diamond}
\def\ox{\mbox{}}
\def\lb{\label}
\def\rel{\;{\rm rel.}\;}
\def\vp{\varepsilon}
\def\ep{\epsilon}
\def\mod{\;{\rm mod}\;}
\def\exp{{\rm exp}\;}
\def\Lie{{\rm Lie}}
\def\diag{{\rm diag}}
\def\im{{\rm im}\;}
\def\Lag{{\rm Lag}}
\def\Gr{{\rm Gr}}
\def\span{{\rm span}}
\def\Stab{{\rm Stab}\;}
\def\sign{{\rm sign}\;}
\def\Supp{{\rm Supp}\;}
\def\Sp{{\rm Sp}\;}
\def\sf{{\rm sf}}
\def\ind{{\rm ind}\;}
\def\rank{{\rm rank}\;}
\def\Sg{{\Sp(2n,\C)}}
\def\even{{\rm even}}
\def\odd{{\rm odd}}

\def\det{{\rm det}\;}
\def\dist{{\rm dist}}
\def\deg{{\rm deg}}
\def\tr{{\rm tr}\;}
\def\ker{{\rm ker}\;}
\def\Vect{{\rm Vect}}
\def\End{{\bf End}}
\def\H{{\bf H}}
\def\k{{\bf k}}
\def\q{{\bf q}}
\def\R{{\bf R}}
\def\C{{\bf C}}
\def\Z{{\bf Z}}
\def\F{{\bf F}}
\def\D{{\bf D}}
\def\A{{\bf A}}
\def\L{{\bf L}}
\def\x{{\bf x}}
\def\y{{\bf y}}
\def\Ga{{\cal G}}
\def\Ha{{\cal H}}
\def\La{{\cal L}}
\def\Pa{{\cal P}}
\def\Ua{{\cal U}}
\def\Va{{\cal V}}
\def\Wa{{\cal W}}
\def\Na{{\cal N}}
\def\E{{\rm E}}
\def\Ja{{\cal J}}
\def\Im{{\rm Im}\;}
\def\m{{\rm m}}
\def\ch{{\rm ch}}
\def\gl{{\rm gl}}
\def\Gl{{\rm Gl}}
\def\Spec{{\bf Spec}\;}
\def\sp{{\rm sp}}
\def\U{{\rm U}}
\def\O{{\rm O}}
\def\F{{\rm F}}
\def\P{{\rm P}}
\def\D{{\rm D}}
\def\T{{\rm T}}
\def\S{{\rm S}}

\title{\bf Holomorphic Equivariant
Cohomology via\\ a Transversal Holomorphic Vector Field
\thanks{Partially supported by the NNSF of China (10271059) and the K. C. Wong Education Foundation}}
\author{Huitao Feng}
\date{$\  $}
\maketitle

\begin{minipage}[t]{14cm}

{\small\noindent{{\bf Abstract}}\hspace{2mm}
In this paper an analytic proof of a generalization of a theorem of Bismut ([Bis1, Theorem 5.1]) is given, which says that,
when $v$ is a transversal holomorphic vector field on a compact complex manifold $X$ with a zero point set $Y$,
the embedding $j:Y\to X$ induces a natural isomorphism between the holomorphic
equivariant cohomology of $X$ via $v$ with coefficients in $\xi$ and the Dolbeault cohomology of $Y$ with coefficients in
$\xi|_Y$, where $\xi\to X$ is a holomorphic vector bundle over $X$.

\noindent{\bf 2000MR Subject Classification 53C, 58J}}
\end{minipage}

\s{Introduction}

\bigskip
For a compact complex manifold $X$, let $T_\C X$ denote the complexification of the real tangent bundle $T_\R X$ of $X$.
Then $T_\C X$ splits canonically as $T_\C X=TX\oplus{\overline{TX}}$, where $TX$ and $\overline{TX}$ are the holomorphic
and the anti-holomorphic tangent bundle of $X$, respectively. Let $\xi\to X$ be a holomorphic vector bundle over $X$.
There exists a natural $\Z$-grading in $\Lambda^\cdot(T_\C^*X)\otimes\xi$ defined by the following decomposition:
$$\Lambda^\cdot(T_\C^*X)\otimes\xi=\bigoplus_{-\dim X\leq r\leq\dim X}\Lambda^{(r)}(T^*_\C X)\otimes\xi,
\eqno(1.1)$$
where
$$\Lambda^{(r)}(T^*_\C X)\otimes\xi=\bigoplus_{q-p=r}\Lambda^q(\overline{T^*X})\otimes(\Lambda^p(T^*X)\otimes\xi).
\eqno(1.2)$$
We will use $\Lambda^{(\cdot)}(T^*_\C X)\otimes\xi$ to denote $\Lambda^\cdot(T^*_\C X)\otimes\xi$ with this $\Z$-grading.
Let $\Omega^{(r)}(X,\xi)$ (resp. $\Omega^{(\cdot)}(X,\xi)$) be the complex vector space of smooth sections of
$\Lambda^{(r)}(T^*_\C X)\otimes\xi$ (resp. $\Lambda^{(\cdot)}(T^*_\C X)\otimes\xi$). Then
$(\Omega^{(\cdot)}(X,\xi)=\bigoplus_r\Omega^{(r)}(X,\xi),\bar\partial^X)$ is a complex and the cohomology
groups $H^{(r)}(X,\xi)$ associated to this complex are direct sums:
$$H^{(r)}(X,\xi)=\bigoplus_{q-p=r}H^{p,q}(X,\xi),
\eqno(1.3)$$
where $H^{p,q}(X,\xi)$ are the usual Dolbeault cohomology groups of $X$ with coefficients in the holomorphic vector bundle
$\xi$.

\bigskip
For any holomorphic vector field $v$ on $X$, set
$$\bar\partial_v^X=\bar\partial^X+i(v):\Omega^{(\cdot)}(X,\xi)\to\Omega^{(\cdot+1)}(X,\xi),
\eqno(1.4)$$
where $i(v)$ is the standard contraction operator defined by $v$. The consideration of operators $\bar\partial^X$ and
$i(v)$ together goes back to [Bot]. Clearly, $(\Omega^{(\cdot)}(X,\xi),\bar\partial_v^X)$
is also a $\Z$-graded complex. Denote the cohomology groups associated to this complex by $H^{(r)}_v(X,\xi)$,
which are called the {\bf holomorphic equivariant cohomology} groups of $X$ via $v$ with coefficients in $\xi$ (cf. [L]).

\bigskip
In [CL], Carrell and Lieberman discussed the relation between Dolbeault cohomology of a connected compact
K$\ddot{\rm a}$hler
manifold $X$ and the zero point set $Y$ of a holomorphic vector field $v$ on $X$ by using Deligne degeneracy criterion and
proved that $H^{(r)}(X,\C)$ vanished for all $|r|>\dim Y$ if $Y\neq\emptyset$. Since
$\dim H^{(\cdot)}_v(X,\C)=\dim H^{(\cdot)}(X,\C)$ in this case (see [L, Theorem 1.3] and [CL]), the corresponding
vanishing results for $H^{(\cdot)}_v(X,\C)$ are also valid. For a general compact complex manifold
$X$ and a transversal holomorphic vector field $v$ on $X$ (see a definition in [Bis1, Sect.5.1]), Liu in [L] proved that
$\dim H^{(\cdot)}(Y,\C)\leq\dim H^{(\cdot)}_v(X,\C)$ by constructing an injective homomorphism
$\alpha_r:H^{(r)}(Y,\C)\to H^{(\cdot)}_v(X,\C)$. Moreover, Liu got a counting formula for $\dim H_v^{(0)}(X,\C)$ in term of
the multiplicities of the zero points of $v$ if the zero points of $v$ are discrete. Under the assumption that the zero
point of $v$ is nodegenerate, motivated by Witten's deformation idea ([W]), he also sketched an analytic proof of his
formula in [L, Sect.7] by examining the behavior of a natural deformation $D^X_T$ of the Riemann-Roch operator on $X$
as $T\to\infty$.

A more general result in this direction is due to Bismut. By using the technique of spectral sequences in
[Bis1, Theorem 5.1], Bismut proves that, when $v$ is a transversal holomorphic vector field on $X$ with a zero point set $Y$,
the embedding $j:Y\to X$ induces naturally a quasi-isomorphism
$j^*:(\Omega^{(\cdot)}(X,\C),{\bar\partial}^X_v)\to(\Omega^{(\cdot)}(Y,\C),{\bar\partial}^Y)$.

\bigskip
In this paper, we will give an analytic proof of the following fairly straightforward generalization of the Bismut's
theorem [Bis1, Theorem 5.1]:

\bigskip
{\bf Theorem 1.1} {\it Let $v$ be a transversal holomorphic vector field on a
compact complex manifold $X$ with a zero point set $Y$. Then
$$j^*:(\Omega^{(\cdot)}(X,\xi),{\bar\partial}^X_v)\to(\Omega^{(\cdot)}(Y,\xi|_Y),{\bar\partial}^Y)
\eqno(1.5)$$
is a quasi-isomorphism.}

\bigskip
Following Witten's deformation idea ([W]) as in [L, Sect.7], we will also work with a deformation $D^X_T$ of the twisted
Riemann-Roch operator on $X$ by $\xi$ but the whole proof now is heavily based on the analytic localization techniques
developed by Bismut and Lebeau (cf. [BL]), since the analysis involved here is much complicated than the situation in
[L, Sect.7]. A key point in our proof is to express $D^X_T$ and to trivialize the bundle on which $D^X_T$ acts by using
the Bismut connection $\nabla^{-B}$ (see [Bis2, II. {\it b)}]). We should point out that if the trivialization is
made by a lifting of the holomorphic Hermitian connection $\nabla^{TX}$ on $TX$, then an extra term coming from the torsion
of $\nabla^{TX}$ will enter to the final operator on $Y$ in the Bismut-Lebeau localization process of $D^X_T$ as
$T\to\infty$. But it is not clear to us that the extra term is zero for a general complex manifold $X$. Consequently, we
could only obtain an equality related to the involved operators at the index level.

\s {A deformed twisted Riemann-Roch operator $D^X_T$ and its local behavior near $Y$}

\bigskip
This section is divided into three parts.
In {\bf a)} we introduce a deformation $D^X_T$ of the twisted Riemann-Roch operator by $\xi$ via $v$, which has been used
in [L, Sect.7] in the case of $\xi=\C$.
In {\bf b)} we recall the definition of the Bismut connection $\nabla^{-B}$ (cf. [Bis2, II. {\it b)}]) and
express $D_T^X$ in this connection by a direct application of [Bis2, Theorem 2.2].
In {\bf c)} we study the local behavior of the deformed operator $D^X_T$ near the submanifold $Y$ following [BL, Sect.8],
in which the Bismut connection will play an essential role.

\bigskip
\noindent{\bf a) A deformed twisted Riemann-Roch operator $D^X_T$}

\bigskip
Let $X$ be a compact complex manifold of $\C$-dimension $n$.
For any $T\in\R$, we consider the following deformed operator
$$\bar\partial^X_T=\bar\partial+Ti(v):\Omega^{(\cdot)}(X,\xi)\to\Omega^{(\cdot+1)}(X,\xi)
\eqno(2.1)$$
and the deformed complex $(\Omega^{(\cdot)}(X,\xi),\bar\partial_T^X)$. One verifies easily as the proof of [L, Lemma 1.1]
that the cohomologies associated to this deformed complex do not depend on $T\neq 0$.

Let $g^{TX}$ (resp. $g^\xi$) be a Hermitian metric on $X$ (resp. $\xi$). By the standard procedure there is an induced
Hermitian metric $\lj\;,\;\rj_{\Lambda^{(\cdot)}(T^*_\C X)\otimes\xi}$ on $\Lambda^{(\cdot)}(T^*_\C X)\otimes\xi$.
Let $d\upsilon_X$ denote the Riemannian volume element of $(X,g^{TX})$. Then for $s_1,s_2\in\Omega^{(\cdot)}(X,\xi)$,
$$\lj\lj s_1,s_2\rj\rj={1\over{(2\pi)^{\dim_\C X}}}\int_X\lj s_1,s_2\rj_{\Lambda^{(\cdot)}(T^*_\C X)\otimes\xi}d\upsilon_X
\eqno(2.2)$$
defines an $L^2$-Hermitian inner product on $\Omega^{(\cdot)}(X,\xi)$. Let $\bar\partial^{X *}$ and $\bar\partial^{X *}_T$
denote the formal adjoint operators of
$\bar\partial^X$ and $\bar\partial^X_T$, respectively. For any $U\in T_\C X$, let $U^*$ be an element in $T^*_\C X$
defined by $g^{TX}(U,\cdot)$. Clearly, $\bar v^*$ is a (1,0)-form on $X$ and $\bar v^*\wedge$ is the dual operator of
$i(v)$. Moreover, we have
$$\bar\partial^{X *}_T=\bar\partial^{X *}+T\bar v^*\wedge:\Omega^{(\cdot)}(X,\xi)\to\Omega^{(\cdot-1)}(X,\xi).
\eqno(2.3)$$
Set
$$D^X={\sqrt 2}(\bar\partial^X+\bar\partial^{X *}),\quad D^X_T={\sqrt 2}(\bar\partial^X_T+\bar\partial^{X *}_T).
\eqno(2.4)$$
Clearly, $D^X$ is the usual twisted Riemann-Roch operator by the holomorphic bundle $\xi$ and $D^X_T$ is a deformation of
$D^X$ and interchanges $\Omega^{({\rm even})}(X,\xi)$ and $\Omega^{({\rm odd})}(X,\xi)$.
From Hodge theory we have the following isomorphisms
$$\ker(D^X_T)^2|_{\Omega^{(r)}(X,\xi)}\cong H^{(r)}_v(X,\xi).
\eqno(2.5)$$

\bigskip
{\bf Lemma 2.1} {\it For any open neighborhood $\Ua$ of $Y$, there exist
constants $a>0$, $b>0$ and $T_0>0$ such that for any $s\in\Omega^{(\cdot)}(X,\xi)$ with $\Supp s\subset X\setminus \Ua$
and any $T\geq T_0$, one has the following estimate for Sobolev norms,
$$\|D^X_Ts\|^2_0\geq a(\|s\|_1^2+(T-b)\|s\|_0^2).
\eqno(2.6)$$}

\bp An easy computation shows that
$$D^{X,2}_T=D^{X,2}+2T^2|v|^2+2T\left((\bar\partial^X\bar v^*)\wedge+i(\bar\partial^X\bar v^*)\right),
\eqno(2.7)$$
where $i(\bar\partial^X\bar v^*)$ denotes the adjoint operator of $(\bar\partial^X\bar v^*)\wedge$. Note that
$((\bar\partial^X\bar v^*)\wedge+i(\bar\partial^X\bar v^*))$ is a zero order operator and $v\neq 0$
on $X\setminus\Ua$, the lemma follows from the well-known Garding's inequality directly.\qe

\bigskip
By Lemma 2.1 and Hodge theory, we can study $H^{(\cdot)}_v(X,\xi)$ through the behavior of the operator $D^X_T$ near $Y$
for large $T$. Also by Lemma 2.1, it is an easy observation that, when $Y=\emptyset$, the cohomology group
$H^{(\cdot)}_v(X,\xi)$ vanishes. In the following we will always assume that $v$ is transversal and $Y\neq\emptyset$.
Note that generally $Y$ consists of some connected components $Y_k$ with different $\C$-dimensions $l_k$.
When no confusion arises, we always drop the subscripts and simply denote them by $Y$ and $l$, respectively.

\bigskip
\noindent{\bf b) An expression of $D_T^X$ in the Bismut connection}

\bigskip
We first recall the definition of the Bismut connection in [Bis2, II. {\it b)}]. For a complex manifold $X$ with a
Hermitian metric $g^{TX}$, let $\nabla^{TX}$ be the holomorphic Hermitian connection on $TX$. Note that $\nabla^{TX}$
induces naturally an Euclidean connection on $T_\R X$ which preserves the complex structure of $T_\R X$. Let $T_X$ denote
the torsion tensor of the connection $\nabla^{TX}$.
Let $B_X$ be the antisymmetrization of the tensor $(U,V,W)\to {1\over 4}\lj T_X(U,V),W\rj$ and let
$S^{-B_X}$ denote the one form with values in antisymmetric elements of $End(T_\R X)$ which is such that
$$\lj S^{-B_X}(U)V,W\rj=-2B_X(U,V,W),
\eqno(2.8)$$
where $U, V, W\in T_\R X$.
Let $\nabla^{L^X}$ be the Levi-Civita connection on $T_\R X$. Set
$$S_X=\nabla^{TX}-\nabla^{L^X}.
\eqno(2.9)$$
Then $S_X$ is also a one form with values in antisymmetric elements of ${\rm End}(T_\R X)$. The important thing here is that
$S^{-B_X}-S_X$ preserves the complex structure of $T_\R X$ (cf. [Bis2, (2.38)]). Now the Bismut connection $\nabla^{-B_X}$
on $T_\R X$ is defined by (cf. [Bis2, (2.37)])
$$\nabla^{-B_X}=\nabla^{TX}+(S^{-B_X}-S_X)=\nabla^{L_X}+S^{-B_X}.
\eqno(2.10)$$
The Bismut connection $\nabla^{-B_X}$ also preserves the complex structure of $T_\R X$
and so induces naturally a unitary connection on $TX$ and a unitary connection on ${\overline{TX}}$, which are still
denoted by $\nabla^{-B_X}$. Note that when $X$ is K$\ddot{\rm a}$hler, the Bismut connection $\nabla^{-B_X}$ coincides with
the holomorphic Hermitian connection $\nabla^{TX}$. There is a unitary connection on
$\Lambda^{\cdot}({\overline{T^*X}})$ lifted canonically from $\nabla^{-B_X}$, which we still call the Bismut connection
and denote by the same notation $\nabla^{-B_X}$.

Let $Y$ be a complex submanifold of $X$. Let $\pi:N\to Y$ be the normal bundle of $Y$ in $X$. We identify $N$ with the
sub-bundle of $TX|_Y$ orthogonal to $TY$ with respect to the restriction metric $g^{TX|_Y}$ on $TX|_Y$ by $g^{TX}$.
So we have the identification of $C^\infty$ bundles $TX|_Y=TY\oplus N$. Let $g^{TY}$ (resp. $g^N$) be the induced metric
on $TY$ (resp. $N$) from $g^{TX|_Y}$. Let $P^{T_\R Y}$, $P^{N_\R}$ be the orthogonal projection operators from $T_\R X|_Y$
onto $T_\R Y$ and $N_\R$ respectively. Let $j$ denote the embedding of $Y$ into $X$. Then $j^*\nabla^{-B_X}$ is a connection
on $T_\R X|_Y$ preserving the metric $g^{TX|_Y}$ and the complex structure of $T_\R X|_Y$. Moreover,
$P^{T_\R Y}(j^*\nabla^{-B_X})P^{T_\R Y}$ is exactly the Bismut connection on $T_\R Y$ associated to the induce metric
$g^{TY}$, and $P^{N_\R}(j^*\nabla^{-B_X})P^{N_\R}$ is a connection on $N_\R$ preserving the metric $g^N$ and the complex
structure of $N_\R$. Set
$$\nabla^{-B_X,\oplus}=\nabla^{-B_Y}\oplus\left(P^{N_\R}(j^*\nabla^{-B_X})P^{N_\R}\right),
\eqno(2.11)$$
$${\bf A}=j^*\nabla^{-B_X}-\nabla^{-B_X,\oplus}.
\eqno(2.12)$$
Clearly, $\nabla^{-B_X,\oplus}$ is also a connection on $T_\R X|_Y$ preserving the metric $g^{TX|_Y}$ and the complex
structure of $T_\R X|_Y$, and ${\bf A}$ is the second fundamental form of the Bismut connection $\nabla^{-B_X}$.
We still use the same notation $\nabla^{-B_X,\oplus}$ to denote its restriction on ${\overline{TX}|_Y}$
as well as its lifting on $\Lambda^{\cdot}(\overline{T^*X}|_Y)$.

\bigskip
Now we return to our situation and express the deformed twisted Riemann-Roch operator $D^X_T$ in the Bismut connection
by applying [Bis2, Theorem 2.2]. To do this we still need a holomorphic Hermitian connection on the bundle
$\Lambda^{\cdot}(T^*X)\otimes\xi$. Since our problem does not depend on the metrics, we can and will choose a special metric
on the holomorphic bundle $\Lambda(T^*X)$ to simplify the analysis.

Let $\L_v:TX|_Y\to TX|_Y$ be the holomorphic Lie homomorphism defined by $\L_v(u)=[v,u]$ for any $u\in TX|_Y$.
Denote $\L_v(TX|_Y)$ by ${\widetilde N}$. Since $v$ is transversal, $TX|_Y$ splits holomorphically into
$TY\oplus{\widetilde N}$ and $\L_v$ induces an isomorphism from $N$ to $\widetilde N$, which we still denote by $\L_v$.
We introduce a new Hermitian metric ${\tilde g}^{\widetilde N}$ on ${\widetilde N}$ by requiring that
$\L_v:N\to{\widetilde N}$ is unitary fiberwisely. Consequently, we get a new Hermitian metric
${\tilde g}^{TX|_Y}=g^{TY}\oplus{\tilde g}^{\widetilde N}$ on $TX|_Y$. We can and we will extend
${\tilde g}^{TX|_Y}$ to a Hermitian metric ${\tilde g}^{TX}$ on $TX$. We will denote $TX$ (resp. $TX|_Y$) with the metric
${\tilde g}^{TX}$ (resp. ${\tilde g}^{TX|_Y}$) by ${\widetilde{TX}}$ (resp. $\widetilde{TX}|_Y$) to distinguish the same
bundle with different metrics. Let $\nabla^{\widetilde{TX}}$ (resp. $\nabla^{\widetilde{TX}|_Y}$, resp. $\nabla^{\widetilde N}$)
be the holomorphic Hermitian connection on ${\widetilde{TX}}$ (resp. ${\widetilde{TX}}|_Y$, resp. ${\widetilde N}$).
We have the following standard fact:
$$j^*\nabla^{\widetilde{TX}}=\nabla^{\widetilde{TX}|_Y}=\nabla^{TY}\oplus\nabla^{\widetilde N}.
\eqno(2.13)$$
We lift the holomorphic Hermitian connection $\nabla^{\widetilde{TX}}$ (resp. $\nabla^{TY}$) to the holomorphic Hermitian
connection $\nabla^{\Lambda^{\cdot}({\widetilde{T^*X}})}$ (resp. $\nabla^{\Lambda^{\cdot}(T^*Y)}$) on
$\Lambda^{\cdot}({\widetilde{T^*X}})$ (resp. $\Lambda^{\cdot}(T^*Y)$).

Let $g^\xi$ be a Hermitian metric on $\xi$ and let $\nabla^\xi$ be the holomorphic Hermitian connection on $\xi$.
Set $\nabla^{\xi|_Y}=j^*\nabla^\xi$, which is the holomorphic Hermitian connection on $\xi|_Y$. So
$\nabla^{{\Lambda^{\cdot}({\widetilde{T^*X}})}\otimes\xi}
=\nabla^{\Lambda^{\cdot}({\widetilde{T^*X}})}\otimes 1+1\otimes\nabla^\xi$
(resp. $\nabla^{\Lambda^{\cdot}(T^*Y)\otimes\xi|_Y}=\nabla^{\Lambda^{\cdot}(T^*Y)}\otimes 1+1\otimes\nabla^{\xi|_Y}$) is
a holomorphic Hermitian connection on $\Lambda^{\cdot}({\widetilde{T^*X}})\otimes\xi$
(resp. $\Lambda^{\cdot}(T^*Y)\otimes\xi|_Y$). Therefore,
$$\nabla^{B,X}=\nabla^{-B_X}\otimes 1+1\otimes\nabla^{{\Lambda^{\cdot}({\widetilde{T^*X}})}\otimes\xi},
\eqno(2.14)$$
$$\nabla^{B,Y}=\nabla^{-B_Y}\otimes 1+1\otimes\nabla^{\Lambda^{\cdot}(T^*Y)\otimes\xi|_Y}
\eqno(2.15)$$
are unitary connections on the Hermitian vector bundle
$\Lambda^{(\cdot)}(T^*_\C X)\otimes\xi$, $\Lambda^{(\cdot)}(T^*_\C Y)\otimes\xi|_Y$, respectively.

\bigskip
For $U\in(TX, g^{TX})$, set
$$c(U)={\sqrt 2}U^*\wedge,\quad c(\bar U)=-{\sqrt 2}i(\bar U);
\eqno(2.16)$$
and for $U\in({\widetilde{TX}},{\tilde g}^{TX})$, set
$${\hat c}(U)=-{\sqrt{-2}}i(U),\quad {\hat c}(\bar U)=-{\sqrt{-2}}{\bar U}^*\wedge.
\eqno(2.17)$$
We extend the map $c$ (resp. $\hat c$) by $\C$ linearity into the Clifford action of $T_\C X$ (resp. ${\widetilde{T_\C X}}$)
on $\Lambda^{\cdot}(\overline{T^*X})$ (resp. $\Lambda^{\cdot}(\widetilde{T^*X})$).

Let $\{e_1,\ldots,e_{2n}\}$ be an orthonormal basis $T_\R X$. Set
$$c(B_X)={1\over 6}\sum^{2n}_{i,j,k=1}B_X(e_i,e_j,e_k)c(e_i)c(e_j)c(e_k).
\eqno(2.18)$$
Note that
$${1\over 4}\sum^{2n}_{i,j,k=1}\lj S^{-B_X}(e_i)e_j,e_k\rj c(e_i)c(e_j)c(e_k)=-3c(B_X).
\eqno(2.19)$$

Now recall the definition (2.10) and apply [Bis2, Theorem 2.2] directly, we obtain the following expressions of $D^Y$
and $D^X_T$ with respect to the orthonormal basis $\{e^\prime_1,\ldots,e^\prime_{2l}\}$ for $T_\R Y$ and
$\{e_1,\ldots,e_{2n}\}$ for $T_\R X$, respectively:
$$D^Y=\sum^{2l}_{i=1}c(e^\prime_i)\nabla^{B,Y}_{e^\prime_i}+2c(B_Y),
\eqno(2.20)$$
$$D_T^X=\sum^{2n}_{i=1}c(e_i)\nabla^{B,X}_{e_i}+2c(B_X)+{\sqrt{-1}}T({\hat c}(v)+{\hat c}(\bar v)).
\eqno(2.21)$$

\bigskip
\noindent{\bf c) The local behavior of the deformed operator $D^X_T$ near $Y$}

\bigskip
For $y\in Y$ and $Z\in N_{\R,y}$, let $t\in\R\to x_t=\exp^X_y(tZ)\in X$ be the geodesic in $X$ with respect to the
Levi-Civita connection $\nabla^{L^X}$, such that $x_0=y$, ${dx}/{dt}|_{t=0}=Z$.
For $\ep>0$, let $B_\ep=\{Z\in N_\R\mid |Z|<\ep\}.$
Since $X$ and $Y$ are compact, there exists an $\ep_0>0$ such that for $0<\ep<\ep_0$, the map
$(y, Z)\in N_\R\to\exp^X_y(Z)\in X$ is a diffeomorphism from $B_\ep$ to a tubular neighborhood $\Ua_\ep$ of $Y$ in $X$.
From now on, we will identify $B_\ep$ with $\Ua_\ep$ and use the notation $x=(y, Z)$ instead of $x=\exp^X_y(Z)$.

We will make use of the trivialization of $(\Lambda^{(\cdot)}(T_\C^*X)\otimes\xi)|_{\Ua_{\ep_0}}$
by the parallel transport of $(\Lambda^{(\cdot)}(T^*_\C X)\otimes\xi)|_Y$ with respect to the connection
$\nabla^{B,X}$ along the geodesic $t\to(y,tZ)$. The key point here is that this trivialization preserves the metric and the
$\Z$-grading since the Bismut connection is a unitary connection and preserves the complex structure of $T_\R X$.
By using the trivialization of
$\Lambda^{(\cdot)}(T_\C^*X)\otimes\xi$ over $\Ua_{\ep_0}$, we can and will make the identification of
$(\Lambda^{(\cdot)}(T_\C^*X)\otimes\xi)|_{\Ua_{\ep_0}}$
with $\pi^*((\Lambda^{(\cdot)}(T_\C^*X)\otimes\xi)|_Y)|_{B_{\ep_0}}$ and so
we can consider $\nabla^{B,X}$ as a unitary connection on the Hermitian vector bundle
$\pi^*((\Lambda^{(\cdot)}(T_\C^*X)\otimes\xi)|_Y)|_{B_{\ep_0}}$ with the obviously induced metric.
Note that there exists another unitary connection $\nabla^{B,X,\oplus}$ on
$\pi^*((\Lambda^{(\cdot)}(T_\C^*X)\otimes\xi)|_Y)$ defined by
$$\nabla^{B,X,\oplus}=\pi^*\left(\nabla^{-B_X,\oplus}\otimes 1
+1\otimes j^*\nabla^{\Lambda^{\cdot}({\widetilde{T^*X}})\otimes\xi}\right).
\eqno(2.22)$$

\bigskip
Let $d\upsilon_Y$ (resp. $d\upsilon_N$) denote the Riemannian volume element of $(Y,g^{TY})$ (resp. the fibres of
$(N,g^N)$). We define a smooth positive function $k(y,Z)$ on $B_{\epsilon_0}$ by the equation
$d\upsilon_X(y,Z)=k(y,Z)d\upsilon_Y(y)d\upsilon_{N_y}(Z)$ and an $L^2$-Hermitian inner product on $\bf E$ by
$$\lj\lj f,g\rj\rj={1\over{(2\pi)^{\dim_\C X}}}\int_Y\int_{N_{\R,y}}\lj f,g\rj(y,Z)d\upsilon_N(Z)d\upsilon_Y(y),
\eqno(2.23)$$
for any $f,g\in{\bf E}$ with compact support, where ${\bf E}$ denotes the set of smooth sections of
$\pi^*((\Lambda^{(\cdot)}(T^*_\C X)\otimes\xi)|_Y)$ on $N_\R$.  Clearly, $k(y)=k(y,0)=1$ on $Y$ and $k(y,Z)$ has a positive
lower bound on $\Ua_{\ep_0/2}$. If $f\in{\bf E}$ has compact support in $B_{\epsilon_0}$, we can identify $f$ with an
element in $\Omega^{(\cdot)}(X,\xi)$ which has compact support in $\Ua_{\ep_0}$.

\bigskip
Let $TN_\R=T^HN_\R\oplus N_\R$ be the splitting of $TN_\R$ induced by the Euclidean connection
$\nabla^{L^N}=P^{N_\R}(j^*\nabla^{L^X})P^{N_\R}$ on $N_\R$, where $T^HN_\R$ denotes the horizontal part of $TN_\R$.
If $U\in T_\R Y$, let $U^H\in T^HN_\R$ denote the horizontal lift of $U$ in $T^HN_\R$, so that $\pi_*U^H=U$. Let
$$\{e_1,\ldots,e_{2l},e_{2l+1},\ldots,e_{2n}\}
\eqno(2.24)$$
be an orthonormal basis of $T_\R X|_Y$ with $\{e_1,\ldots,e_{2l}\}$ an orthonormal basis of $T_\R Y$ and
$\{e_{2l+1},\ldots,e_{2n}\}$ an orthonormal basis of $N_\R$.

\bigskip
{\bf Definition 2.2} {\it Let $D^H$, $D^N$ be the operators acting on $\bf E$
$$D^H=\sum^{2l}_{i=1}c(e_i)\nabla^{B,X,\oplus}_{e^H_i}+2c(B_Y),\quad
D^N=\sum^{2n}_{\alpha=2l+1}c(e_\alpha)\nabla^{B,X,\oplus}_{e_\alpha}.
\eqno(2.25)$$}

Clearly, $D^N$ acts along the fibres $N_{\R,y}$ as the operator
${\sqrt 2}(\bar\partial^{N_y}+\bar\partial^{N_y*})$.
Note that $D^H$, $D^N$ are self-adjoint with respect to the Hermitian inner product (2.23).

\bigskip
Now we turn to Taylor expansions of $v$ near $Y$ along the geodesic $(y,tZ)$ for $y\in Y$ and $Z\in N_{\R,y}$. Let
$\{w_1,\ldots, w_l, w_{l+1},\ldots, w_n\}$ be a unitary basis for $TX|_Y$ and let $(z^{l+1},\ldots,z^{n})$ denote the
associated holomorphic coordinate system on $N_y$ with $w_\alpha={\sqrt 2}{\partial/{\partial z^\alpha}}$ for
$l+1\leq\alpha\leq n$.
Note that $\L_v:N\to{\widetilde N}$ is unitary fiberwisely. Set
$${\tilde w}_\alpha=\L_v(w_\alpha),\quad{\tilde{\bar w}}_\alpha=\L_v({\bar w}_\alpha),\quad l+1\leq\alpha\leq n.
\eqno(2.26)$$
Hence, $\{w_1,\ldots, w_l, {\tilde w}_{l+1},\ldots, {\tilde w}_n\}$
is a unitary basis for $\widetilde{TX}|_Y$. We use ${\tilde w}^\tau$ to denote the parallel transport of ${\tilde w}$
with respect to the holomorphic Hermitian connection $\nabla^{\widetilde{TX}}$ along the geodesic $(y,tZ)$.
We write $v$ on $\Ua_\ep$ as
$$v(y,Z)={1\over{\sqrt 2}}\left(\sum^{l}_{i=1}v^i(y,Z)w^\tau_i+\sum^n_{\alpha=l+1}v^\alpha(y,Z){\tilde w}^\tau_\alpha\right)
\eqno(2.27)$$
for some smooth functions $v^i$ and $v^\alpha$. Set
$$v_{Y,1}={1\over{\sqrt 2}}\sum^{l}_{i=1}\sum^{n}_{\alpha=l+1}{\partial v^i\over\partial z^\alpha}(y)z^\alpha w_i,
\;v_{N,1}={1\over{\sqrt 2}}\sum^{n}_{\alpha,\beta=l+1}{\partial v^\alpha\over\partial z^\beta}(y)z^\beta{\tilde w}_\alpha,
\eqno(2.28)$$
$$v_{Y,2}={1\over 2{\sqrt 2}}\sum^{l}_{i=1}\sum^{n}_{\alpha,\beta=l+1}
{{\partial}^2v^i\over{{\partial z^\alpha}{\partial z^\beta}}}(y)z^\alpha z^\beta w_i,
\;v_{N,2}={1\over 2{\sqrt 2}}\sum^{n}_{\alpha,\beta,\gamma=l+1}
{{\partial}^2v^\alpha\over{{\partial z^\beta}{\partial z^\gamma}}}(y)z^\beta z^\gamma{\tilde w}_\alpha.
\eqno(2.29)$$
Since $v$ is transversal, we get by the definition (2.26),
$$v_{Y,1}=0,\quad v_{N,1}=-{1\over{\sqrt 2}}\sum^n_{\alpha=l+1}z^\alpha{\tilde w}_\alpha,
\eqno(2.30)$$
and so
$$v(y,Z)=v_{N,1}(y,Z)+v_{Y,2}(y,Z)+v_{N,2}(y,Z)+O(|Z|^3).
\eqno(2.31)$$
Define
$$D^N_T=D^N+{\sqrt{-1}}T({\hat c}\left(v_{N,1})+{\hat c}({\bar v}_{N,1})\right).
\eqno(2.32)$$
A direct and easy computation shows that
$$(D^N_T)^2=-4\sum^n_{\alpha=l+1}{\partial^2\over{{\partial z^\alpha}{\partial{\bar z}^\alpha}}}
+T^2|Z|^2-{\sqrt{-1}}T\sum^n_{\alpha=l+1}(c({\bar w}_\alpha){\hat c}({\tilde w}_\alpha)
+c(w_\alpha){\hat c}({\bar{\tilde w}}_\alpha)).
\eqno(2.33)$$
Set
$$\theta=\sum^n_{\alpha=l+1}w_\alpha^*\wedge{\bar{\tilde w}}_\alpha^*.
\eqno(2.34)$$
Clearly, $\theta$ is a well-defined smooth section of $\Lambda^{(\cdot)}({\bar N}^*\otimes{\widetilde N}^*)$ over $Y$ of the
degree 0. Now we have the following analogue of [BL, Proposition 7.3]:

\bigskip
{\bf Lemma 2.3} {\it Take $T>0$. Then for any $y\in Y$, the operator $(D^N_T)^2$ acting on
$\Gamma\left(\pi^*\Lambda^{(\cdot)}({\bar N^*_y}\oplus{\widetilde N_y})\right)$ over $N_y$ is nonnegative with
the kernel $\C\{\beta_y\}$, where
$$\beta_y=\exp(\theta_y-{T\over 2}|Z|^2),\quad |\exp\theta_y|_{\Lambda^{(\cdot)}({\bar N^*_y}\oplus{\widetilde N_y})}
=2^{{(n-l)}/2}.
\eqno(2.35)$$
Moreover, the nonzero eigenvalues of $(D^N_T)^2$ are all $\geq TA$ for some positive constant $A$ which can be chosen to be
independent of $y$.}

\bp The proof of the lemma is standard (cf. [BL, Sect.7, (7.10)--(7.13)]; also cf. [Z1, Chapter 4, Sect. 4.5]).
\qe

\bigskip
For any $y\in Y$, $Z\in N_{\R,y}$, let $\tau U$ denote the parallel transport of $U\in T_{\R,y}X$ with respect to the
Levi-Civita connection $\nabla^{L^X}$ along the geodesic $(y,tZ)$. Note that we have identified the bundle
$\pi^*((\Lambda^{(\cdot)}(T_\C^*X)\otimes\xi)|_Y)|_{B_{\ep_0}}$ with the bundle
$(\Lambda^{(\cdot)}(T_\C^*X)\otimes\xi)|_{\Ua_{\ep_0}}$ by trivializing the later bundle along the geodesic $(y,tZ)$
by using the connection $\nabla^{B,X}$. The Clifford action of $c((\tau U)(y,tZ))$ on $\Lambda^{(\cdot)}(T_\C^*X)\otimes\xi$
is generally not constant along the geodesic $(y,tZ)$. This is different from the situation in [BL, Sect.8], where the
connection on the related bundle is the lifting of the Levi-Civita connection on $TX$ since the manifold $X$ is
K$\ddot{\rm a}$hler. Hence, to obtain an analogue of [BL, Theorem 8.18], we need to work out the difference between
$c((\tau U)(y,Z))$ and the constant Clifford action $c(U)$ on
$\pi^*((\Lambda^{(\cdot)}(T_\C^*X)\otimes\xi)|_Y)|_{B_{\ep_0}}$.
Since $\nabla^{-B_X}=\nabla^{L^X}+S^{-B_X}$ is unitary, we know that
$$[\nabla^{B,X}_Z,c((\tau U)(y,tZ))]_{|_{t=0}}=c(\nabla^{-B_X}_Z(\tau U)(y,tZ))_{|_{t=0}}=c(S^{-B_X}(Z)U),
\eqno(2.36)$$
thus
$$c((\tau U)(y,Z))=c(U)+c(S^{-B_X}(Z)U)+O(|Z|^2).
\eqno(2.37)$$

\bigskip
Set with respect to the basis (2.24):
$$M={1\over 2}\sum^{2l}_{i,j}\sum^{2n}_{\alpha=2l+1}\lj{\bf A}(e_i)e_j,e_\alpha\rj c(e_i)c(e_j)c(e_\alpha)
-{1\over 2}\sum^{2n}_{\alpha=2l+1}(e_\alpha k)c(e_\alpha),
\eqno(2.38)$$
$$c(B^\prime(y))={1\over 2}\sum^{2l}_{i,j=1}\sum^{2n}_{\alpha=2l+1}B_X(e_i,e_j,e_\alpha)c(e_i)c(e_j)c(e_\alpha)$$
$$+{1\over 6}\sum^{2n}_{\alpha,\beta,\gamma=2l+1}B_X(e_\alpha,e_\beta,e_\gamma)c(e_\alpha)c(e_\beta)c(e_\gamma),
\eqno(2.39)$$
$$c({B^\prime}^\prime(y))={1\over 2}\sum^{2l}_{i=1}\sum^{2n}_{\alpha,\beta=2l+1}
B_X(e_i,e_\alpha,e_\beta)c(e_i)c(e_\alpha)c(e_\beta).
\eqno(2.40)$$
One verifies easily that
$$c(B_X(y))=c(B_Y(y))+c(B^\prime(y))+c({B^\prime}^\prime(y)).
\eqno(2.41)$$

\bigskip
Now we have the following analogue of [BL, Theorem 8.18], which describes the local behavior of $D^X_T$ near $Y$.
Comparing to [BL, Theorem 8.18, (8.58)], some new terms enter into the following theorem.

\bigskip
{\bf Theorem 2.4} {\it As $T\to +\infty$, then
$$k^{1/2}D^X_Tk^{-1/2}=D^H+D^N_T+M_T+c({B^\prime}^\prime)+T{\sqrt{-1}}{\hat c}(v_{Y,2}+{\bar v}_{Y,2})+{\bf S}+R_T,
\eqno(2.42)$$
where
$$M_T=M+c(B^\prime)+T{\sqrt{-1}}{\hat c}(v_{N,2}+{\bar v}_{N,2}),
\eqno(2.43)$$
$${\bf S}=-\sum^{2n}_{i=1}c(e_i)\nabla^{B,X,\oplus}_{P^{N_\R}S^{-B_X}(Z)e_i},
\eqno(2.44)$$
$$R_T=O(|Z|\partial^H+|Z|^2\partial^N+|Z|+T|Z|^3),
\eqno(2.45)$$
and $\partial^H$, $\partial^N$ represent horizontal and vertical differential operators, respectively.}

\bp Let $\{\tau e_1,\ldots,\tau e_{2n}\}$ be the parallel transport of the basis (2.24) with respect to the Levi-Civita
connection $\nabla^{L^X}$ along the geodesic $(y,tZ)$ for $y\in Y$ and $Z\in N_{\R,y}$.
From (2.21), we have
$$D_T^X=\sum^{2n}_{i=1}c(\tau e_i)\nabla^{B,X}_{\tau e_i}+2c(B_X)+{\sqrt{-1}}T({\hat c}(v)+{\hat c}(\bar v)).
\eqno(2.46)$$
We identify $\overline{TX}$ with $\pi^*(\overline{TX}|_Y)$ over $\Ua_{\ep_0}$ by trivializing $\overline{TX}$
with respect to the Bismut connection $\nabla^{-B_X}$ along the geodesic $(y,tZ)$ and set
$$\Gamma=\nabla^{-B_X}-\nabla^{-B_X,\oplus}.
\eqno(2.47)$$
Let $\Gamma^\wedge$ denote the lifting action of $\Gamma$ on $\pi^*(\Lambda^{\cdot}(\overline{T^*X})|_Y)$.
For any $y\in Y$, we find by (2.12)
$$\sum^{2n}_{i=1}c(e_i)\Gamma_y^\wedge(e_i)
={1\over 2}\sum^{2l}_{i,j}\sum^{2n}_{\alpha=2l+1}\lj{\bf A}_y(e_i)e_j,e_\alpha\rj c(e_i)c(e_j)c(e_\alpha).
\eqno(2.48)$$
Furthermore, recall (2.13) and then we get
\begin{eqnarray*}
k^{1/2}D_T^Xk^{-1/2}&=&\sum^{2n}_{i=1}c(\tau e_i)\nabla^{B,X,\oplus}_{\tau e_i}+2c(B_X)+
{\sqrt{-1}}T({\hat c}(v)+{\hat c}(\bar v))\\
&&+\sum^{2n}_{i=1}c(e_i)\Gamma_y^\wedge(e_i)-{1\over 2}\sum^{2n}_{\alpha=2l+1}(e_\alpha k)(y)c(e_\alpha)+O(|Z|).
\end{eqnarray*}
Note that (2.31), (2.38), (2.41), (2.43) and (2.48), we have
\begin{eqnarray*}
k^{1/2}D_T^Xk^{-1/2}&=&\sum^{2n}_{i=1}c(\tau e_i)\nabla^{B,X,\oplus}_{\tau e_i}
+2c(B_Y)+T{\sqrt{-1}}{\hat c}(v_{N,1}+{\bar v}_{N,1})\\
&&+M_T+c({B^\prime}^\prime)+T{\sqrt{-1}}{\hat c}(v_{Y,2}+{\bar v}_{Y,2})+O(|Z|+T|Z|^3).
\end{eqnarray*}
By (2.37) and the expansion of $\tau e_i$ along $(y,tZ)$ in the proof of [BL, Theorem 8.18], especially [BL, (8.80), (8.84)],
we have
\begin{eqnarray*}
\sum^{2n}_{i=1}c(\tau e_i)\nabla^{B,X,\oplus}_{\tau e_i}&=&\sum^{2l}_{i=1}c(e_i)\nabla^{B,X,\oplus}_{e_i^H}
+\sum^{2n}_{\alpha=2l+1}c(e_\alpha)\nabla^{B,X,\oplus}_{e_\alpha}\\
&&+\sum^{2n}_{i=1}c(S^{-B_X}(Z)e_i)\nabla^{B,X,\oplus}_{e_i}+O(|Z|\partial^H+|Z|^2\partial^N),
\end{eqnarray*}
and then by the definition of $D^H$ and $D^N_T$,
\begin{eqnarray*}
k^{1/2}D_T^Xk^{-1/2}&=&D^H+D^N_T+M_T+c({B^\prime}^\prime)+T{\sqrt{-1}}{\hat c}(v_{Y,2}+{\bar v}_{Y,2})\\
&&+\sum^{2n}_{i=1}c(S^{-B_X}(Z)e_i)\nabla^{B,X,\oplus}_{e_i}+O(|Z|\partial^H+|Z|^2\partial^N+|Z|+T|Z|^3).
\end{eqnarray*}
But
\begin{eqnarray*}
\sum^{2n}_{i=1}c(S^{-B_X}(Z)e_i)\nabla^{B,X,\oplus}_{e_i}&=&-\sum^{2n}_{i=1}c(e_i)\nabla^{B,X,\oplus}_{S^{-B_X}(Z)e_i}\\
&=&-\sum^{2n}_{i=1}c(e_i)\nabla^{B,X,\oplus}_{P^{N_\R}S^{-B_X}(Z)e_i}
-\sum^{2n}_{i=1}c(e_i)\nabla^{B,X,\oplus}_{P^{T_\R Y}S^{-B_X}(Z)e_i}\\
&=&{\bf S}+O(|Z|\partial^H),
\end{eqnarray*}
from which we complete the proof of the theorem.
\qe

\s {The proof of Theorem 1.1}

In this section, we prove Theorem 1.1 by using Bismut-Lebeau's techniques in [BL, Sect.9, Sect.10, a)].

\bigskip
For any $\mu\geq 0,$ let $\E^{\mu}$ (resp. ${\bf E}^{\mu}$, resp. $\F^{\mu}$) be the set of sections of
$\Lambda^{(\cdot)}(T_\C^*X)\otimes\xi$ on $X$ (resp. of $\pi^*((\Lambda^{(\cdot)}(T_\C^*X)\otimes\xi)|_Y)$ on the
total space of $N$, resp. of $\Lambda^{(\cdot)}(T_\C^*Y)\otimes\xi|_Y$ on Y) which lie in the $\mu$-th Sobolev spaces.
Let $\|\quad\|_{\E^{\mu}}$ (resp. $\|\quad\|_{{\bf E}^{\mu}}$, resp. $\|\quad\|_{\F^{\mu}}$) be the Sobolev norm on
$\E^{\mu}$ (resp. ${\bf E}^{\mu}$, resp. $\F^{\mu}$). We always assume that the norms $\|\quad\|_{\E^0}$ (resp.
$\|\quad\|_{{\bf E}^0}$, resp. $\|\quad\|_{\F^0}$ is the norm associated with the scalar products on the corresponding
bundles).

Let $\gamma:\R\to [0, 1]$ be a smooth even function with $\gamma(a)=1$ if $|a|\leq {1\over 2}$ and $\gamma(a)=0$
if $|a|\geq 1.$ For any $y\in Y$, $Z\in N_y$ and $\ep\in (0,\ep_0)$, where $\ep_0$ is chosen as in Section 2, ${\bf c)}$, set
$$\gamma_\ep(Z)=\gamma\left({{|Z|}\over\ep}\right),\quad
\alpha_T=\int_{N_{\R,y}}\gamma^2_\ep(Z)\exp(-T|Z|^2){{d\upsilon_{N_y}(Z)}\over{(2\pi)^{\dim_\C N}}}.
\eqno(3.1)$$
Clearly, $\alpha_T$ does not depend on $y\in Y$ and $\alpha_T=O({1\over{T^{n-l}}})$.

For $\mu\geq 0,T>0,$ define linear maps $I_T: {\rm F}^{\mu}\to{\bf E}^{\mu}$ and $J_T:{\rm F}^{\mu}\to\E^\mu$ by
$$I_T u=\left({1\over{2^{n-l}\alpha_T}}\right)^{1/2}\gamma_\ep(Z)(\pi^*u)\beta_y ,\; J_T u=k^{-1/2}I_T u,
\quad\forall\;u\in {\rm F}^\mu.
\eqno(3.2)$$
It is easy to see that $I_T$, $J_T$ are isometries from ${\rm F}^0$ onto their images.
For $\mu\geq 0, T>0,$ let ${\bf E}^{\mu}_T$ (resp. $\E^\mu_T$) be the image of ${\rm F}^{\mu}$ in ${\bf E}^{\mu}$
(resp. $\E^\mu$) under $I_T$ (resp. $J_T$) and let ${\bf E}^{0,\perp}_T$ (resp. ${\E}^{0,\perp}_T$) be the orthogonal
complement of ${\bf E}^0_T$ (resp. ${\E}^0_T$) in ${\bf E}^0$ (resp. ${\E}^0$) and let $p_T$, $p_T^{\perp}$
(resp. $\bar p_T, \bar p_T^\perp$) be the orthogonal projection operators from ${\bf E}^{0}$ (resp. $\E^0$)
onto ${\bf E}^{0}_T, {\bf E}^{0,\perp}_T$ (resp. $\E^0_T, \E^{0,\perp}_T$), respectively. Set
$$\E^{\mu,\perp}=\E^{\mu}\cap\E_T^{0,\perp}.
\eqno(3.3)$$
Then $\E^0$ splits orthogonally into
$$\E^0=\E^0_T\oplus\E^{0,\perp}_T.
\eqno(3.4)$$
Since the map
$s\in {\bf E}^{0}\to k^{-1/2}s \in \E^0  $ is an isometry, we see that the map $s \to k^{-1/2}s$ identifies the Hilbert
space ${\bf E}^0_T$ and $\E^0_T.$ Corresponding to the decomposition (3.4) we set:
$$D_{T,1}=\bar p_TD^X_T\bar p_T,\quad D_{T,2}=\bar p_TD^X_T\bar p_T^\perp,
\quad D_{T,3}=\bar p_T^\perp D^X_T\bar p_T,\quad D_{T,4}=\bar p_T^\perp D^X_T\bar p_T^\perp.
\eqno(3.5)$$
Then
$$D^X_T = D_{T,1} + D_{T,2} + D_{T,3} + D_{T,4}.
\eqno(3.6)$$

\bigskip
We have the following analogue of [BL, Theorem 9.8].

\bigskip
{\bf Lemma 3.1} {\it The following formula holds on
$\Gamma(\Lambda^{(\cdot)}(T_\C^*Y)\otimes\xi|_Y)$ as $T\to +\infty$
$$J_T^{-1}D_{T,1}J_T=D^Y +O({1\over\sqrt T}),
\eqno(3.7)$$
where $O({1\over\sqrt T})$ is a first order differential operator with smooth coefficients dominated by $C/\sqrt T$.}

\bp Note that the action of the operator $M_T$ on $\pi^*((\Lambda^{(\cdot)}(T^*_\C X)\otimes\xi)|_Y)$
interchanges $\pi^*(\Lambda^{(\rm even)}({\bar N}^*\oplus{\widetilde N}^*))$ and
$\pi^*(\Lambda^{(\rm odd)}({\bar N}^*\oplus{\widetilde N}^*))$, we get
$$p_TM_Tp_T=0.
\eqno(3.8)$$
Note that $B_X$ is antisymmetric and $\lj c(e_\alpha)c(e_\beta)\beta_y,\beta_y\rj=0$ for any $\alpha,\beta$
with $2l+1\leq\alpha<\beta\leq 2n$, we obtain
$$p_Tc({B^\prime}^\prime)p_T=0.
\eqno(3.9)$$
Since $B_X$ is antisymmetric, by (2.8) we get $\lj S^{-B_X}(Z)e_i,Z\rj=0$ and so
$$\nabla^{B,X,\oplus}_{P^{N_\R}S^{-B_X}(Z)e_i}\left(\gamma_\ep(Z)\exp(\theta_y-{{T|Z|^2}\over 2})\right)
=\exp(\theta_y-{{T|Z|^2}\over 2})(P^{N_\R}S^{-B_X}(Z)e_i)\gamma_\ep(Z).
\eqno(3.10)$$
From the equality above we can prove easily the following estimate for some uniformly positive constant $C$ and
any $s\in{\bf E}^1$:
$$\|p_T{\bf S}s\|_{{\bf E}^0}\leq{C\over{\sqrt T}}\|s\|_{{\bf E}^1}.
\eqno(3.11)$$
Since $\int_{\C}e^{-T|z|^2}z^2dzd{\bar z}=0$ and $(\gamma^2_\ep-1)$ vanishes on a symmetric domain containing $0$,
we have for any $u\in {\rm F}$ that
$$I^{-1}_Tp_T(\hat c(v_{Y,2})+{\hat c}({\bar v}_{Y,2}))
\left({1\over{2^{n-l}\alpha_T}}\right)^{1/2}\gamma_\ep(Z)(\pi^*u)\beta_y=O({1\over{T^{3/2}}}).
\eqno(3.12)$$
On the other hand, note that $\beta_y$ is of constant length on $Y$, we get for $1\leq i\leq 2l$ that
$\lj\nabla^{B,X,\oplus}_{e_i}\beta,\beta\rj=0$ and so
$$I_T^{-1}p_TD^Hp_TI_T=D^Y.
\eqno(3.13)$$
One can then proceed as in [BL, Proof of Theorem 9.8] and use (3.8)--(3.12) to complete the proof of Lemma 3.1
easily.
\qe

\bigskip
Note that the estimate (3.12) and proceed as the proof of Theorem 9.10, Theorem 9.11 and Theorem 9.14 in
[BL, Sect.9], one can prove the following lemma without any new difficulty.

\bigskip
{\bf Lemma 3.2} {\it There exist $C_1>0$, $C_2>0$ and $T_0>0$ such that for any $T\geq T_0$,
$s\in\E_T^{1,\perp}$ and $s^\prime\in\E^1_T,$ then
$$\|D_{T,2}s\|_{\E^0}\leq C_1\left({\|s\|_{\E^1}\over\sqrt T}+\|s\|_{\E^0}\right),
\eqno(3.14)$$
$$\|D_{T,3}s^\prime\|_{\E^0}\leq C_1\left({\|s^\prime\|_{\E^1}\over\sqrt T}+\|s^\prime\|_{\E^0}\right),
\eqno(3.15)$$
$$\|D_{T, 4}s\|_{\E^0}\geq C_2(\|s\|_{\E^1}+\sqrt T\|s\|_{\E^0}).
\eqno(3.16)$$}

\bigskip
Let $\Spec(D^Y)$ denote the spectrum of $D^Y$. Choose $c>0$ such that $\Spec(D^Y)\cap [-2c, 2c]\subset\{0\}.$
Let $\delta=\{\lambda\in\C: |\lambda|=c\}$. Let $\E_c(T)$ denote the direct sum of the eigenspaces of $D^X_T$ with
eigenvalues lying in $[-c, c]$. Then $\E_c(T)$ is a finite dimensional subspace of $\E^0.$ Let $P_{T,c}$ denote the
orthogonal projection from $J_T(\ker(D^Y))$ to $\E_c(T)$. By Lemma 3.1 and Lemma 3.2, we have the following analogue of
[BL, (9.156)] (also see [TZ, Proposition 4.4] for a proof without using the norm in [BL, Sect.9, Definition 9.17] and the
distance in [BL, Sect.9, Definition 9.22]):

\bigskip
{\bf Theorem 3.3} {\it There exist $c>0$ and $T_0>0$ such that for any $T\geq T_0$, the projection
$$P_{T,c}:J_T(\ker(D^Y))\to \E_c(T)
\eqno(3.17)$$
is an isomorphism.}

\bigskip
Now to prove Theorem 1.1 we only need to prove that when $T$ large enough, $D^X_T$ has no nonzero small eigenvalues or
equivalently, to prove the following equality:
$$E_c(T)=\ker(D^X_T).
\eqno(3.18)$$

\bigskip
Let $Q$ denote the orthogonal projection from $\Omega^{(\cdot)}(Y,\xi)$ to $\ker(D^Y)$. Then we have the following analogue
of [Z2, Theorem 1.10] (also see [BL, Theorem 10.1, (10.4)]):

\bigskip
{\bf Theorem 3.4} {\it There exist $c>0$, $C>0$, $T_1>0$ such that for any $T\geq T_1$, any $\sigma\in\ker(D^Y)$,
$$\|(2^{n-l}\alpha_T)^{1/2}Qj^*P_{T,c}J_T\sigma-\sigma\|_0\leq {C\over{\sqrt T}}\|\sigma\|_0.
\eqno(3.19)$$}

\bp The proof of [Z2, Theorem 1.10], which is a modified version of the proof of [BL, Theorem 10.1, (10.4)], is carried out
here with the identity [Z2,(1.34)] in the proof of the [Z2, Theorem 1.10] replaced by the following equality
$$j^*{1\over{2\pi{\sqrt{-1}}}}\int_\delta k^{-{1/2}}\gamma_\ep
{{(\pi^*\sigma)\beta_y}\over\lambda}d\lambda=\sigma.
\eqno(3.20)$$
The proof of (3.20) is similar to that of the identity [BL, (10.29)].\qe

\bigskip
Note that $j^*\beta_y=1$ is crucial in the proof of (3.19). It is no longer true for the case of the analytic proof of Morse
inequalities of Witten ([W]) since in that case the contribution of the bundle $\Lambda(N^*)$ to the kernel of $D^N_T$ is a
pure $p$-form around each critical point of index $p$ (cf. [Z1, Chapter 5, 6]) and its pull-back by $j^*$ vanishes on $Y$.
Consequently, $j^*$ can not be a quasi-isomorphism at all in that case.

\bigskip
{\bf Proof of Theorem 1.1.} The proof of Theorem 1.1 now is similar to that in [Z2, Sect.1, {\bf e)}]. Note that
the trick used in Zhang's proof and so ours is inspired by Braverman ([Br, Sect.3]). First of all, we know
from Theorem 3.3 that
$$P_{T,c}J_T:\ker(D^Y)\to\E_c(T)
\eqno(3.21)$$
is an isomorphism when $T$ is very large.
Take $\alpha\in\E_c(T)$. Then ${\bar\partial}^X_T\alpha\in\E_c(T)$. By the above discussion, there exists
$\sigma\in\ker(D^Y)$ such that
$${\bar\partial}^X_T\alpha=(2^{n-l}\alpha_T)^{1/2}P_{T,c}J_T\sigma.
\eqno(3.22)$$
From (3.22) and that $j^*$ is a quasi-homomorphism, i.e. $j^*{\bar\partial}_T^X={\bar\partial}^Yj^*$, we have
$$(2^{n-l}\alpha_T)^{1/2}Qj^*P_{T,c}J_T\sigma
=Qj^*{\bar\partial}^X_T\alpha=Q{\bar\partial}^Yj^*\alpha=0.
\eqno(3.23)$$
From (3.23) and (3.19), we get
$$\|\sigma\|_0\leq{C\over{\sqrt T}}\|\sigma\|_0,
\eqno(3.24)$$
and so $\sigma=0$ as $T$ large enough.
Thus, when $T$ is large enough, we have that
$${\bar\partial}^X_T|_{\E_c(T)}=0.
\eqno(3.25)$$
From (3.25) and Theorem 3.3, we have that when $T$ is large enough,
$$\dim\ker(D^X_T)=\dim\E_c(T)=\dim\ker(D^Y).
\eqno(3.26)$$
Now by Theorem 3.4,
$$j^*:\ker(D^X_T)\to\ker(D^Y)
\eqno(3.27)$$
is clearly an injective and so an isomorphism from (3.26).
\qe

\bigskip

\bigskip
\noindent{{\bf Acknowledgements} The author would like to thank Professors Jean-Michel Bismut and Xiaonan Ma
for many helpful discussions from which the author benefits a lot. This work was done while the author was visiting the
Institut des Hautes $\acute{\rm E}$tudes Scientifiques in Bures-Sur-Yvette. He would like to thank Professor
Jean-Pierre Bourguignon and the IHES for their hospitality and support.

\bigskip

\bigskip
\centerline{{\bf References}}

\bigskip
\noindent[Bis1] J. -M. Bismut, Holomorphic and de Rham torsions, {\it preprint}.

\noindent[Bis2] J. -M. Bismut, A local index theorem for non K$\ddot{\rm a}$hler manifolds,
{\it Math. Ann.} 284(1989), 681-699.

\noindent[BL] J. -M. Bismut and G. Lebeau, Complex immersions and Quillen metrics,
{\it Publ. Math. IHES.} V.74 (1991).

\noindent[Bot] R. Bott, A residue formula for holomorphic vector fields, {\it J. Diff. Geom.} 1(1967), 311-330.

\noindent[Br] M. Braverman, Cohomology of the Mumford quotient, {\it Progress in Math.} 198(2001), 47-59.

\noindent[CL] J. B. Carrell and D. I. Lieberman, Holomorphic vector fields and Kaehler manifolds,
{\it Invent. Math.} 21 (1973), 303-309.

\noindent[L] K. Liu, Holomorphic equivariant cohomolgy, {\it Math. Ann.} 303(1995), 125-148.

\noindent[TZ] Y. Tian and W. Zhang, Quantization formula for symplectic
manifolds with boundary, {\it Geom. Funct. Anal.} 9(1999), 596-640.

\noindent[W] E. Witten, Supersymmetry and Morse theory, {\it J. Diff. Geom.} 17(1982), 661-692.

\noindent[Z1] W. Zhang, Lectures on Chern-Weil Theory and Witten Deformations, {\it Nankai Tracts in Mathematics}, V.4,
World Scientific, Singapore, 2001.

\noindent[Z2] W. Zhang, A holomorphic quantization formula in singular reduction, {\it Commun. Contemp. Math.} 1(1999),
No.3, 281-293.

\bigskip

\bigskip
Feng: College of Mathematical Sciences, Nankai University, Tianjin, 300071, China;\\
fht@nankai.edu.cn

Current address:

Feng: Institut Des Hautes $\acute{\rm E}$tudes Scientifiques, Bures-Sur-Yvette, 91440, France;\\
feng@ihes.fr}

\end{document}